\documentclass{article}[12pt]
\usepackage{amsmath,amssymb,amsthm}
\advance\topmargin by -1in
 
\def\Bbb#1{{\mathbb #1}}
\def\Cal#1{{\cal #1}}
\def\Frak#1{\mathfrak{#1}}
 
\def\then{\;\Longrightarrow\;}
\def\iff{\;\Longleftrightarrow\;}

\def\Span{\mathop{\rm Span}\nolimits}

\def\lc{\mathop{\rm lc}\nolimits}

\def\Span{\mathop{\rm Span}\nolimits}

\def\ord{\mathop{\rm ord}\nolimits}

\def\lc{\mathop{\rm lc}\nolimits}
\def\ind{\mathop{\rm ind}\nolimits}
\def\supp{\mathop{\rm supp}\nolimits}

\newcommand{\GR}{Gr\"{o}bner}

\def\Forall{\mbox{\ for each }}
\def\Where{\mbox{\ where }}

\def\And{\mbox{\ and }}

\def\Jm{\mathop{\rm Im}\nolimits}
\def\qed{\ifmmode\squareforqed\else{\unskip\nobreak\hfil
\penalty50\hskip1em\null\nobreak\hfil\squareforqed
\parfillskip=0pt\finalhyphendemerits=0\endgraf}\fi}
\def\squareforqed{\hbox{\rlap{$\sqcap$}$\sqcup$}}

\theoremstyle{plain} 
\newtheorem{Theorem}{Theorem}
\newtheorem{Corollary}[Theorem]{Corollary}

\theoremstyle{definition} 
\newtheorem{Definition}[Theorem]{Definition}

\theoremstyle{remark} 
\newtheorem{Construction}[Theorem]{Construction}

\newtheorem{Remark}[Theorem]{Remark}

\newcommand{\FF}{\mathbb{F}}
\newcommand{\cL}{\mathcal{L}}

\renewcommand\labelitemi{\labelitemii}
 
\title{Macaulay, Lazard and the Syndrome Variety}
\author{M. Ceria\thanks{The author expresses her heartful thanks to Teo Mora for the many fruitful discussions and suggestions on this topic.}}

\begin{document}

\maketitle

\begin{abstract}
In this paper we consider the four syndrom varieties ${\sf Z}_e^\times$, i.e. the set of all error locations corresponding to errors of weight $w, 0\leq w\leq 2$,
 ${\sf Z}_{ns}^\times$ , the set of all {\em non spurious} error locations corresponding to errors of weight $w, 0\leq w\leq 2$,  ${\sf Z}_+^\times $, the set of all non-spurious  error locations  corresponding to errors of weight $w, 1\leq w\leq 2$,  ${\sf Z}_2^\times $, the set of all non-spurious  error locations  corresponding to errors of weight $w= 2$, associated to 
an up-to-two errors correcting binary cyclic codes. Denoting $J_\ast:=\mathcal{I}({\sf Z}_\ast)$, the ideal of these syndrome varieties,   ${\sf N}_\ast := {\bf N}(J_\ast)$ the \GR\ escalier of $J_\ast$ w.r.t. the lex ordering with  $x_1<x_2<z_1<z_2$,
 $\Phi_\ast : {\sf Z}_\ast \to {\sf N}_\ast$ a Cerlienco-Mureddu correspondence, and $G_*$ a minimal Groebner basis of the ideal $J_\ast$, the aim of the paper is, assuming to know 
  the structure of the order ideal ${\sf N}_2$ and  a Cerlienco Mureddu Correspondence 
to deduce with elementary arguments ${\sf N}_\ast$, $G_\ast$ and $\Phi_\ast$ for $\ast\in\{e,ns,+\}$.

The tools are Macaulay's trick and Lazard's formulation of Cerlienco-Mureddu correspondence.
\end{abstract}

\section{Introduction}

Classical, decoding of BCH codes $C\subset{\FF}_q^n$ are based on solving the {\em key equation} \cite{Berl} $$\sigma(x)S(x)\equiv \omega(x) \pmod{x^{2t}},$$
where 
\begin{itemize}
\item $\alpha\in \FF_q[\alpha]=:{\FF}_{q^m}$ is a primitive $n^{th}$-root of the unity;
\item $S(x) = \sum_{i=1}^{2t}s_{i}x^{i-1}$, $s_i := \sum_{j=1}^{\mu} e_{\ell_j}\alpha^{i\ell_j}$, is the syndrome polynomial associated to the error $\sum_{j=1}^{\mu}
e_{\ell_j}\alpha^{\ell_j}, \mu\leq t$;
\item $\sigma(x)=\prod_{j=1}^{\mu} (1-x\alpha^{\ell_j})$ is the classical error locator polynomial\footnote{The recent mood, not depending on the network technology of the Sixties, and thus  not needing the key equation,
prefers consider the plain error locator polynomial $\prod_{j=1}^{\mu} (x-\alpha^{\ell_j})$.};
\item and $\omega(x)=\sum_{j=1}^{\mu}e_{\ell_j}\alpha^{\ell_j}\prod\limits_{\genfrac{}{}{0pt}{}{i=1}{i\neq j}}^\mu (1-x\alpha^i)$ is Newton's {\em error evaluator polynomial}
\end{itemize}

In 1990 Cooper \cite{l7, l8} suggested to use \GR \ basis
computation in order to decode cyclic codes: let $C$ be a binary
BCH code correcting up to $t$ errors,
$\bar{s}=(s_1,\dots,s_{2t-1})$ be the syndrome vector associated
to a received word.  Cooper's idea consisted in interpreting the
error locations $z_1,\ldots z_t$ of $C$ as the roots of the syndrome equation
system:
$$f_i := \sum_{j=1}^t z_j^{2i-1} - s_{2i-1} = 0, \,\,\,1\leq i \leq t,$$
and, consequently, the plain error locator polynomial  as the monic generator
$g(z_1)$ of the principal ideal
$$\left\{\sum_{i=1}^{t} g_i f_i, g_i \in \FF_{2}(s_1,\ldots,s_{2t-1})[z_1,\ldots,z_t]\right\}\bigcap \FF_{2}(s_1,\ldots,s_{2t-1})[z_1],$$
which was computed via the elimination property of lexicographical \GR \ bases.

In a series of papers \cite{CRHT2, CRHT3, CRHT} Chen et al.
improved and generalized Cooper's approach to decoding. In
particular, for a $q$-ary $[n,k,d]$ cyclic code, with correction
capability $t$, they made the following alternative proposals:
\begin{enumerate}
\item denoting, for an error with weight $\mu$,
 $z_1,\ldots,z_\mu$ the error locations, $y_1,\ldots,y_\mu$ the error values,
$s_1,\ldots,s_{n-k}\in \FF_{q^m}$ the associated syndromes,
they interpreted \cite{CRHT2}
 the coefficients of the plain error locator polynomial
 as the elementary symmetric functions $\sigma_j(z_1,\ldots,z_\mu)$
 and the syndromes as the {\em Waring functions},  $s_i = \sum_{j=1}^\mu
 y_jz_j^{i}$,
and suggested to deduce the $\sigma_j$'s from the (known) $s_i$'s
via a \GR\ basis computation of the ideal generated by the Newton
identities; a similar idea was later developed in \cite{faugere,faugere2}.
 \item They considered \cite{CRHT3} the {\em syndrome variety} (Definition~\ref{SynVar})
 {\footnotesize{$$\left\{(s_1,\ldots,s_{n-k}, y_1,\ldots,y_t, z_1,\ldots,z_t)\in (\FF_{q^m})^{n-k+2t} : s_i = \sum_{j=1}^\mu
y_jz_j^{i}, \,1 \leq i \leq n-k \right\}$$}} and proposed to
deduce via a Gr\"obner basis pre-computation in
$$\FF_{q}[x_1,\ldots,x_{n-k}, y_1,\ldots,y_t, z_1,\ldots,z_t]$$ a series of polynomials
$g_\mu(x_1,\ldots,x_{n-k},Z), \mu \leq t$ such that, for any error
with weight $\mu$ and associated syndromes $s_1,\ldots,s_{n-k}\in
\FF_{q^m}$, $g_\mu(s_1,\ldots,s_{n-k},Z)$ in $\FF_{q^m}[Z]$ is the
plain error locator polynomial.

This approach was improved in a series of paper \cite{Caboara,LY} culminating with \cite{locator} which, specializing Gianni-Kalkbrener Theorem \cite{gianni,kalk}, stated Theorem~\ref{OSALA1} below.
\end{enumerate}
For a survey of this {\em Cooper Philosophy} see \cite{Cooper} and on Sala-Orsini locator \cite{COST}.

\medskip\noindent
Recently the same problem has been reconsidered within the frame of {\em Grobner-free Solving} \cite{Lu, Mou, Lu, FGLML}, explicitly expressed and sponsored in the book \cite[Vol.3,40.12,41.15]{SPES}; such approach aims to avoid the computation of a \GR\ bases of a (0-dimensional) ideal $J\subset{\Cal P}$ in favour of  combinatorial algorithms, describing instead the structure of the algebra ${\Cal P}/J$. In particular, given a finite set of distinct points $\mathbf{X}=\{P_1,...,P_N\}$ and, for each point $P_i$ the related primary ${\sf q}_i$ described via suitable functionals  $\{\ell_{i1},\ldots,\ell_{ir_i}\}$, the aim is to describe the combinatorial structure of
both the ideal ${\sf I} = \cup_i {\sf q}_i$ and the {\em algebra} $\FF[{\bf N}({\sf I})]={\Cal P} \setminus {\sf I}$
avoiding Buchberger Algorithm and Groebner bases of ${\sf I}$ and even Buchberger reduction modulo ${\sf I}$ in favour of combinatorial algorithms, as variations of Cerlienco-Mureddu Algorithm \cite{Ce,CeMuBC,CeMu,CeMu2,CeMu3,FRR} and Lundqvist interpolation formula \cite{Lu}, dealing with $\FF[{\bf N}({\sf I})]$ and $\{\ell_{ij}, 1\leq i\leq N, 1\leq j\leq r_i\}$.

In particular  \cite{AcaHelp}, given the syndrome variety 
$${\sf Z}=\left\{(c+d,c^l+d^l,c,d), c,d\in\FF_{2^m}^\ast,c\neq d\right\}$$ of a 
binary cyclic code  
$C$ over $\FF_{2^m}$, with length $n=2^m-1$ and primary defining set $S_C=\{1,l\}$, and denoted ${\Cal I}({\sf Z})\subset\FF_2[x_1,x_2,z_1,z_2]$ the ideal of points of ${\sf Z}$, produced with good complexity, via Cerlienco-Mureddu Algorithm \cite{CeMu,CeMu2,CeMu3},  
the set ${\bf N} :={\bf N}({\Cal I}({\sf Z}))$ and proved that
the related \GR\ basis has the shape
$$G=(x_1^n-1,g_2,z_2+z_1+x_1,g_4)$$
where (see \cite{locator})
$$g_2=\frac{x_2^{\frac{n+1}{2}}-x_1^\frac{n+1}{2}}{x_2-x_1}=
x_2^{\frac{n-1}{2}}+\sum_{i=1}^{\frac{n-1}{2}} \binom{\frac{n-1}{2}}{i}x_1^ix_2^{\frac{n-1}{2}-i}$$
and $g_4=z_1^2-\sum_{t\in{\bf N}} c_tt$ is Sala-Orsini general error locator polynomial. 

Such result allowed \cite{AcaHelp} to remark (applying Marinari-Mora Theorem \cite{Iran,Mother2,Ce0}) that, for decoding, it is sufficient to 
express ${\sf Z}$ as
$${\sf Z}=\left\{(c+ca^{i},c^l+c^la^{li},c,ca^{i}),c\in\FF_{2^m}^\ast, 1\leq i\leq n-1 \right\}$$
and 
compute the  
polynomial --
{\em half error locator polynomial} (HELP) --
$$h(x_1,x_2,z_1):=z_1-\sum_{t\in{\bf H}} c_tt \Where {\bf H}:=\{x_1^{i}x_2^{j}, 0\leq i< n, 0\leq j< \frac{n-1}{2}\}
$$
which satisfies 
$$h(c(1+a^{2j+1}),c^l(1+a^{l(2j+1)}),z_1)=z_1-c, \Forall c\in\FF_{2^m}^\ast, 0\leq j< \frac{n-1}{2},$$ 
the other error $ca^{2j+1}$ been computable via the polynomial $z_2+z_1+x_1\in G$ as
$z_2:=x_1-z_1=(c+ca^{2j+1})-c=ca^{2j+1}.$

Such polynomial can be easily obtained with good complexity via Lundqvist interpolation formula \cite{Lu} on the set of points
$$\{(c+ca^{2j+1},c^l+c^la^{l(2j+1)},c), c\in\FF_{2^m}^\ast, 0\leq j< \frac{n-1}{2}\}.$$
Experiments showed that in that setting HELP has a very sparse formula, which has been proved in \cite{AcaHelp}:
$$h(x_1,x_2,z_1)=z_1+\sum_{i=1}^{\frac{n-1}{2}} a_i x_1^{(n+1-li) \mod n}x_2^{(i-1)\mod \frac{n-1}{2}}$$
where  the unknown coefficients can be deduced by Lundqvist interpolation on the set of points
$$\{(1+a^{2j+1},1+a^{l(2j+1)},1), 0\leq j< \frac{n-1}{2}\}$$
and on the monomials $\{x_1^{(n+1-li) \mod n}x_2^{(i-1)\mod \frac{n-1}{2}},1\leq i< \frac{n+1}{2}\}$,
thus granting an ${\Frak O}(n^{\bf ?})$ combinatorial pre-processing produces a linear decoder.
\medskip\noindent

This suggested to perform the same investigation to the general case $n\mid 2^m-1$.
Denoting by
\begin{itemize}
 \item  
 $a$ a primitive $(2^m-1)^{\rm th}$ root of unity so that ${\mathbb F}_{2^m}={\mathbb Z}_2[a]$, $\alpha:=\frac{2^m-1}{n}$ and 
\item $b:=a^{\alpha}$ a primitive $n^{\rm th}$ root of unity, 
\item ${\Cal R}_n:=\{e\in{\mathbb F}_{2^m} : e^n=1\}$   the set of the $n^{\rm th}$ roots of unity and
\item ${\Cal S}_n :={\Cal R}_n\sqcup\{0\}$;
\end{itemize}
let us consider the following sets of points
\begin{description}
\item[] ${\sf Z}_2 :=\{(c+d,c^l+d^l,c,d), c,d\in{\Cal R}_n, c\neq d\}$,$\#{\sf Z}_2^\times=n^2-n$;
\item[] ${\sf Z}_+ :=\{(c+d,c^l+d^l,c,d), c,d\in{\Cal S}_n, c\neq d\}$,$\#{\sf Z}_+^\times=n^2+n$,
\item[] ${\sf Z}_{ns} :=\{(c+d,c^l+d^l,c,d), c,d\in {\Cal S}_n\}\setminus\{(0,0,c,c),c\in{\Cal R}_n\}$,$\#{\sf Z}_{ns}^\times=n^2+n+1$,
\item[] ${\sf Z}_e :=\{(c+d,c^l+d^l,c,d), c,d\in {\Cal S}_n\}$,$\#{\sf Z}_e^\times=(n+1)^2$,
\end{description}
and let us denote, for  $\ast\in\{e,ns,+,2\}$, 
\begin{itemize}
 \item $J_\ast:=\mathcal{I}({\sf Z}_\ast)$, 
 \item ${\sf N}_\ast := {\bf N}(J_\ast)$ the \GR\ escalier of $J_\ast$ w.r.t. the lex ordering with  $x_1<x_2<z_1<z_2$ and 
 \item $\Phi_\ast : {\sf Z}_\ast \to {\sf N}_\ast$ a Cerlienco-Mureddu correspondence.
\end{itemize}

The aim of the paper is, assuming to know 
\begin{enumerate}
\renewcommand\theenumi{{\rm (\alph{enumi})}}  
\item the structure of the order ideal ${\sf N}_2$, $\#{\sf N}_2=n^2-n$, i.e. a minimal basis $\{t_1,\ldots,t_r\}, t_i:=x_1^{a_i}x_2^{b_i},$ of the monomial ideal 
${\Cal T}\setminus{\sf N}_2={\bf T}({\Frak I}({\sf Z}_2))$,

\item a Cerlienco Mureddu Correspondence $\Phi_2: {\sf N}_2 \to {\sf Z}_2$
\end{enumerate}
to deduce with elementary arguments ${\sf N}_\ast$ and $\Phi_\ast$ for $\ast\in\{e,ns,+\}$\footnote{and, since it is completely obvious to deduce it, also the irrelevant \GR\ bases $G_\ast$ of the ideal $J_\ast$.}.

\section{Notations}\label{Notations}
\noindent 
${\mathbb F}$ denotes an arbitrary field, $\overline{\mathbb F}$ denotes
its algebraic closure and ${\mathbb F}_q$ denotes a finite field of size
$q$ (so $q$ is implicitly understood to be a power of a prime) and
${\Cal P} :=  {\mathbb F}[X] := {\mathbb F}[x_1,\ldots,x_n]$  the polynomial ring over the field ${\mathbb F}$.

Let
${\Cal T}$ be the set of terms in ${\Cal P}$, {\em id est}
$${\Cal T} := \{x_1^{a_1}\cdots x_n^{a_n} :
(a_1,\ldots,a_n)\in {\Bbb N}^{n}\}.$$

If $t=x_1^{\gamma_1}\cdots x_n^{\gamma_n}\in{\Cal T}$, then $\deg(t)=\sum_{i=1}^n
\gamma_i$ is the \emph{degree} of $t$ and, for each $h\in \{1,...,n\}$, $\deg_h(t):=\deg_{x_h}(t):=\gamma_h$ is the $h$-\emph{degree} of $t$.

A \emph{semigroup ordering} $<$ on $\mathcal{T}$  is  a total ordering
 such that  $$t_1<t_2 \Rightarrow st_1<st_2,\, \Forall s,t_1,t_2
\in \mathcal{T}.$$ For each semigroup ordering $<$ on $\mathcal{T}$,  we can represent a polynomial
$f\in \mathcal{P}$ as a linear combination of terms arranged w.r.t. $<$, with
coefficients in the base field ${\mathbb F}$:
$$f=\sum_{t \in \mathcal{T}}c_tt=\sum_{t \in \mathcal{T}}c(f,t)t=\sum_{i=1}^s c(f,t_i)t_i:\,
c(f,t_i)\in
{\mathbb F}\setminus \{0\},\, t_i\in \mathcal{T},\, t_1>...>t_s.$$  
For each such $f$ 
its {\em support} is
$\supp(f)  := \{\tau\in {\Cal T}  : c(f,t) \neq 0\},$
its {\em leading term} is 
the term ${\bf T}_<(f) := \max_<(\supp(f))=t_1$, its {\em
leading coefficient} is $\lc_<(f) := c(f,t_1)$ and its {\em leading monomial}
is ${\bf M}_<(f) := \lc_<(f){\bf T}_<(f)=c(f,t_1)t_1.$ 
When $<$ is understood we will drop the subscript, as in ${\bf T}(f) = {\bf T}_<(f)$.

A \emph{term ordering} is a semigroup ordering  such that  $1$ is lower 
than every variable or, equivalently, it is a \emph{well ordering}.\\
Given an ordered set of varibles $[x_1,\ldots,x_n]$ we consider the \emph{lexicographical ordering} 
induced
by  $x_1<...<x_n$, i.e:
$$ x_1^{\gamma_1}\cdots x_n^{\gamma_n}<_{Lex} x_1^{\delta_1}\cdots
x_n^{\delta_n} \Leftrightarrow \exists j\, \vert  \,
\gamma_j<\delta_j,\,\gamma_i=\delta_i,\, \forall i>j, $$
which is a term ordering. Since we do not consider any 
term ordering other than Lex, we drop the subscript and denote it by $<$ 
instead of $<_{Lex}$.

The assignement of a finite set of terms 
$${\sf G} := \{\tau_1,\ldots,\tau_\nu\}\subset{\Cal T}, \tau_i = x_1^{a^{(i)}_1}\cdots x_n^{a^{(i)}_n}$$
defines  a partition ${\cal T} = {\sf T}\sqcup {\sf N}$ of ${\Cal T}$ in two parts:
\begin{itemize}
\item ${\sf T} := \{\tau \tau_i : \tau\in {\Cal T}, 1\leq i \leq\nu \}$
which is a {\em semigroup ideal}, 
{\em id est} a subset ${\sf T}\subset{\Cal T}$  such that 
$$\tau\in {\sf T}, t\in{\cal T}\then t\tau\in {\sf T};$$
\item the {\em normal set} ${\sf N}:=  {\Cal T} \setminus {\sf T}$
which is an {\em order ideal},
{\em id est} a subset ${\sf N}\subset{\Cal T}$ such that 
$$\tau\in {\sf N},  t\in{\cal T}, t\mid\tau\then t\in {\sf N},$$
\end{itemize}

For any set $F \subset {\Cal P}$, write
\begin{itemize}
\item ${\bf T}\{F\} := \{{\bf T}(f) : f\in F\};$
\item ${\bf M}\{F\} := \{{\bf M}(f) : f\in F\}$;
\item ${\bf T}(F) := \{\tau {\bf T}(f) : \tau\in {\Cal T}, f\in F\},$ a semigroup ideal;
\item ${\bf N}(F) :={\Cal T} \setminus {\bf T}(F),$ an order ideal;
\item ${\Bbb I}(F) = \langle F \rangle$ the ideal generated by $F$.
\item ${\mathbb F}[{\bf N}(F)] :=\Span_{\mathbb F}({\bf N}(F)).$
\end{itemize}

Given an ideal $J\subset{\Cal P}$, denote ${\sf G}(J)$ the
minimal set of generators of ${\bf T}(J)$;
we denote by ${\sf T}(J):={\bf T}(J)$
the semigroup ideal and ${\sf{N}}(J):={\bf N}(J)$ the order ideal introduced by the partition
${\cal T} = {\bf T}(J)\sqcup {\bf N}(J)= {\sf T}(J)\sqcup {\sf N}(J)$
defined by ${\sf G}(J)$; ${\sf{N}}(J)$ will be called the \emph{Gr\"obner escalier} of $J$.

\medskip

\noindent Let $\mathbf{X}=\{P_1,...,P_N\} \subset {\mathbb F}^n$ be a finite set of distinct points,
$$P_i:=(a_{1,i},...,a_{n,i}),\, i=1,...,N.$$
We call
$${\Cal I}(\mathbf{X}):=\{f \in \mathcal{P}:\, f(P_i)=0,\, \forall i\},$$
the \emph{ideal of points} of $\mathbf{X}$.\\
If we are interested in the \emph{ordered set}, instead of its support $\mathbf{X}$, 
we denote it by  $\underline{\mathbf{X}}=[P_1,...,P_N]$.

For any (0-dimensional, radical) ideal $J\subset{\Cal P}$ and any extension field $E$ of ${\mathbb F}$, let
${\Cal V}_E(J)$ be the (finite) rational points of $J$ over $E$. We also write
${\Cal V}(J) = {\Cal V}_{\overline{\mathbb F}}(J)$. We have the obvious duality between ${\Cal I}$ and
${\Cal V} = {\Cal V}_{\overline{\mathbb F}}.$

\begin{Definition} An ordered finite set $\underline{\mathbf{X}}=[P_1,...,P_N] \subset {\mathbb F}^n$ of points and an ordered finite set
$\underline{Q}:[q_1,\cdots,q_N]\subset\mathcal{P}$ of polynomials are said to be
\begin{itemize}
\item triangular iff $q_i(P_j)=\begin{cases}0&i< j,\cr 1&i=j;\cr\end{cases}$ 
\item biorthogonal iff $q_i(P_j)=\begin{cases}0&i\neq j,\cr 1&i=j.\cr\end{cases}$
\end{itemize}
\end{Definition}

\begin{Definition}\label{AA}  For an ideal ${\sf I}\subset{\Cal P}$, a finite set $G\subset{\sf I}$ will be called a {\em Gr\"obner
basis} of ${\sf I}$ if 
${\bf T}(G) = {\bf T}({\sf I}),$ that is,
${\bf T}\{G\} := \{{\bf T}(g) : g\in G\}$ generates  
${\bf T}({\sf I}) = {\bf T}\{{\sf I}\}$. 
\end{Definition}

\section{Ingredients}
\subsection{Cerlienco--Mureddu Correspondence}
Given an {\em ordered} finite set of points $\underline{\mathbf{X}}=[P_1,...,P_N] \subset {\mathbb F}^n$, set ${\sf N}:= {\bf N}({\Cal I}(\mathbf{X}))$ and, for each $\sigma, 1\leq \sigma \leq N$,  $\underline{\mathbf{X}}_\sigma :=[P_1,\ldots,P_\sigma]$.

A Cerlienco--Mureddu Correspondence \cite{Iran, Mother2} for $\underline{\mathbf{X}}$ is a bijection $\Phi : \underline{\mathbf{X}} \to {\sf N}$ such that, for each $\sigma, 1\leq \sigma \leq N$,
$\underline{\sf N}_\sigma :[\Phi(P_1),\ldots,\Phi(P_\sigma)]$ is the Gr\"obner escalier of the ideal ${\bf N}({\Cal I}(\mathbf{X}_\sigma)).$

Such correspondence is the  output of a combinatorial algorithm proposed by Cerlienco and Mureddu \cite{CeMu,CeMu2,CeMu3} and produces exactly the same values ${\bf N}({\Cal I}(\mathbf{X}_\sigma))$ as M\"oller Algorithm \cite{MMM2}. 

Cerlienco--Mureddu Algorithm is {\em inductive} and thus has complexity $\mathcal{O}\left(n^2N^2\right)$ (see \cite{FRR}), but it has the advantage of being {\em iterative}, in the sense that, given an ordered set of points 
$\underline{\mathbf{X}}=[P_1,...,P_N]$, its related escalier 
$\underline{\sf N} = {\bf N}({\Cal I}(\mathbf{X}))$
and correspondence $\Phi_{\underline{\mathbf{X}}} : \underline{\mathbf{X}} \to \underline{\sf N}$, for any point $Q\notin\underline{\mathbf{X}}$ it returns a term 
$\tau\in\mathcal{T}$ such that, denoting $\underline{\mathbf{Y}}$ the ordered set $\underline{\mathbf{Y}}:=[P_1,...,P_N,Q]$, 
\begin{itemize}
 \item $ {\bf N}({\Cal I}(\mathbf{Y}))= {\bf N}({\Cal I}(\mathbf{X}))\sqcup\{\tau\}$,
 \item $\Phi_{\underline{\mathbf{Y}}}(P_i)=\Phi_{\underline{\mathbf{X}}}(P_i)$ for all $i$ and
 $\tau=\Phi_{\underline{\mathbf{Y}}}(Q)$.
\end{itemize}

In order to produce an {\em iterative} procedure without paying the inductivity complexity of \cite{CeMu}, \cite{CeMuBC} applies the Bar Code \cite{Ce, CeJ}
which describes in a compact way the combinatorial strucure of a (non necessarily 0-dimensio\-nal) 
ideal; the Bar Code actually 
allows to remember and reed those data which Cer\-lien\-co-Mured\-du algorithm is forced to inductively recompute. Actually
the application of the Bar Code allows to compute the lexicographical Gr\"obner escaliers ${\sf N}_\sigma$ and the related 
Cer\-lien\-co-Mured\-du correspondences, with iterative complexity $\mathcal{O}(N\cdot N\log(N)n)$.

\subsection{Macaualay's trick}
In \cite[p.458]{Mac} Macaualay proposed the following
\begin{Construction}
Given a finite set of generators ${\sf G} :=\{t_1,\ldots,t_p\}$ of a 0-dimensional monomial ideal
${\sf T}\subset {\Cal T}\subset {\mathcal P},$
where, setting $\delta :=\max\{\deg(t) : t\in{\sf G}\}$ and assuming that the field ${\mathbb F}$ contains a copy of $\mathbb{Z}_\delta$,
associate
\begin{itemize}
 \item to each term $t_\iota := x_1^{e_{\iota1}}
\cdots x_n^{e_{\iota n}}, 1\leq \iota\leq p$, the polynomial
 $p_\iota := \prod_{l=1}^n \prod_{j=0}^{e_{\iota l}-1} x_l-j$,
  \item to each term $n_i := x_1^{e_{i1}}
\cdots x_n^{e_{in}}, 1\leq i\leq m =:\#{\sf N}$ the point $P_i := \left(e_{i1},
\cdots ,e_{in}\right) \in \mathbb{Z}_\delta^n$.
\end{itemize}
\end{Construction}
Macaulay's aim was to show that each function $H: \mathbb{N}\to\mathbb{N}$, which was satisfying some precise bounds stated by him, was the Hilbert function of a (0-dimensional, radical) ideal and in particular the vanishing ideal of a (finite) set of separate points; and he remarked that, given such a function $H$, if ${\sf G}$ satisfies
$$H(d) = \#\{n\in {\sf N},\deg(n)=d\},\Forall d\leq \delta \And H(d)=0 \Forall d>\delta$$ then
\begin{enumerate}
 \item $\mathbb{I}\left(p_\iota : 1\leq \iota \leq p\right) = {\Cal I}\left(P_i : 1\leq i\leq m\right)$.
\end{enumerate}

Our interest is toward the further remarks \cite{NG2} that
\begin{enumerate}\setcounter{enumi}{1}
\item the map 
$$\Phi : {\sf N} \to \mathbf{X} := \left\{P_i : 1\leq i\leq m\right\} : n_i \mapsto P_i\Forall i, 1\leq i\leq m$$
is a Cerlienco-Mureddu Correspondence;
\item $\left\{p_\iota : 1\leq \iota \leq p\right\}$ is a Gr\"obner basis of the ideal it generates w.r.t. each termordering;
\item we remark moreover that, if $\underline{\sf N}$ (and $\underline{\mathbf{X}}$) is enumerated \cite[II,p.553]{SPES} in such a way that, for each $j, 1\leq j\leq m$, the set 
$\underline{\sf N}_j := \left[n_i : 1\leq i\leq j\right]$ is an order ideal and we associate 
\begin{itemize}
  \item to each term $n_i := x_1^{e_{i1}}
\cdots x_n^{e_{in}}$ the polynomial
 $q_i := \prod_{l=1}^n \prod_{j=0}^{e_{il}-1} \frac{x_l-j}{e_{il}-j}$,
\end{itemize}
than each set $Q_j :=\left[q_i : 1\leq i\leq j\right]$ is a triangular set for each set $\underline{\mathbf{X}}_j := \left[P_i : 1\leq i\leq j\right]$. 
\end{enumerate}

\subsection{Lazard's Cerlienco-Mureddu Correspondence}

While many algorithms  are available  which, given a finite set of points $\underline{\mathbf{X}}$ produce the normal set ${\sf N}\subset {\Cal T}$ s.t. 
${\sf N} = {\bf N}(\mathcal{I}({\mathbf{X}})$ \cite{BM,CeMu,CeMu2,CeMu3,MMM2,FRR,CeMuBC} and a Cerlienco--Mureddu Correspondence $\Phi : \underline{\sf N} \to \underline{\mathbf{X}}$ \cite{CeMu,CeMu2,CeMu3,MMM2, CeMuBC}, since we are interested in applying it to the elementary case $n=2$, we make reference to the first and stronger instance proposed by Lazard \cite{L}:
denoting $\pi : {\mathbb F}^2 \to {\mathbb F}$ the projection s.t. $\pi(a_1,a_2) = a_1$ for each $(a_1,a_2)\in {\mathbb F}^2$, given a finite set of distinct points
$\mathbf{X}\subset {\mathbb F}^2$ set
\begin{itemize}
\renewcommand\labelitemi{\labelitemii}
 \item $\{a_0,\ldots,a_{r-1}\} := \pi({\mathbf{X}}),$
 \item $d(i) := \#\left\{P\in{\mathbf{X}} : \pi(P)=a_i \right\}$;
\end{itemize}
after renumerating the $a_i$s we can assume $d(0)\geq d(1) \geq \ldots \geq d(r-1)$; thus there are values $b_{ij}, 0\leq i < r, 0\leq j < d(i),$ such that
$${\mathbf{X}}=\left\{(a_i,b_{ij}), 0\leq i < r, 0\leq j < d(i)\right\};$$
then
\begin{itemize}
 \item ${\bf N}(\mathbb{I}({\mathbf{X}})=\left\{X_1^iX_2^j:  0\leq i < r, 0\leq j < d(i)\right\}$
 \item $\Phi(a_i,b_{ij}) = X_1^iX_2^j$ is a Cerlienco-Mureddu Correspondence.
\end{itemize}

\subsection{Syndrome Variety, spurious roots and syndrome map}
We begin briefly recalling the standard notation on cyclic codes, needed to understand what follows, making reference to \cite{Berl,PH}

Let $C$ be a $[n,k,d]_q$ a $q$-ary cyclic code with length $n$, dimension $k$ and distance $d$. We denote by $g(x)\in \FF_q[x]$ its {\em generator 
polynomial}, remarking that $\deg(g)=n-k$ and $g \mid x^n-1$. Let $\FF_{q^m}$ be the splitting field of $x^n-1$ over $\FF_q$.\\
If $a$ is a primitive $n$-th root of unity, the {\em complete defining set} of $C$ is 
$$S_C=\{j\vert g(a^j)=0,\, 0 \leq j \leq n-1\}.$$
This set  is completely partitioned in cyclotomic classes, so 
we can pick an element for each such class, getting a set $S\subset S_C$, uniquely identifying the code. This set $S$ is a {\em primary defining set } of $C$.\\
If $H$ is a parity-check matrix of $C$, $\mathbf{c}$ is a codeword (i.e. $\mathbf{c} \in C$),  
$\mathbf{e}\in (\FF_q)^n$  an error vector and $\mathbf{v}=\mathbf{c}+ \mathbf{e}$ a received vector, the vector $\mathbf{s}\in 
(\FF_{q^m})^{n-k}$ such that its transpose  $\mathbf{s^T}$ is $\mathbf{s^T}=H\mathbf{v^T}$ is called {\em syndrome vector}.
We call {\em correctable syndrome} a syndrome vector corresponding to an error of weight $\mu \leq t$, where $t$ is the {\em error correction capability} of the code, i.e. the maximal number of errors that the code can correct.

The notion of {\em syndrome variety} (see \cite{Cooper}) was formalized in \cite{CRHT3} in their approach to decoding $q$-ary $[n,k,d]$ cyclic codes, with correction capability $t$.
\begin{Definition}\label{SynVar} For such a cyclic code, the {\em syndrome variety} is the set of points
$$\mathbf{V}:=\left\{(s_1,\ldots,s_{n-k}, y_1,\ldots,y_t, z_1,\ldots,z_t)\in (\FF_{q^m})^{n-k+2t} : s_l = \sum_{j=1}^\mu
y_jz_j^{l}, \,1 \leq i \leq n-k \right\}$$
where for an error $(s_1,\ldots,s_{n-k}, y_1,\ldots,y_t, z_1,\ldots,z_t)\in \mathbf{V}$ with weight $\mu\leq t$ and
$$y_{\mu+1}=\cdots=y_t=0, \quad z_{\mu+1}=\cdots=z_t=0,$$ 
$z_1,\ldots,z_\mu$ represent the {\em  error locations}, $y_1,\ldots,y_\mu$ the {\em error values},
$s_1,\ldots,s_{n-k}\in \FF_{q^m}$ {\em the associated syndromes}.
\end{Definition}

\begin{Definition} For such a cyclic code, and $\mu\leq t$ the {\em plain error locator polynomial} is the polynomial $\prod_{j=1}^\mu (X-z_i)$
\end{Definition}

\begin{Definition} For such a cyclic code, we call {\em syndrome map}, the function
$$\phi :  (\FF_{q^m})^t\to (\FF_{q^m})^{n-k} : (z_1,\ldots,z_t)\mapsto \left(\sum_{j=1}^\mu
y_jz_j,\ldots,\sum_{j=1}^\mu
y_jz_j^{l}\right).$$
\end{Definition}

\begin{Definition} \cite{CRHT3,locator} A point $(s_1,\ldots,s_{n-k}, y_1,\ldots,y_t, z_1,\ldots,z_t)\in\mathbf{V}$ is said {\em spurious} if there 
are at least two values $z_i,z_j, 1\leq i \neq j \leq \mu$ such that $z_i=z_j\neq0.$
\end{Definition}

\subsection{GELP}

In 2005, Orsini and Sala \label{OSL} proved that each cyclic code $C$ possesses a polynomial $\cL$ of the form
$$\cL=z^t+\sum_{i=1}^ta_{t-i}(x_1,...,x_{n-k})z^{t-i}\in \FF_q[x_1,...,x_{n-k}][z],$$
which
satisfies the following condition
 \begin{quote}
given
a syndrome vector ${s}=({s}_1,\ldots, {s}_{n-k}) \in ({\Bbb
F}_{q^m})^{n-k}$ corresponding to an error with weight $\mu \leq
t$, then its $t$ roots are the $\mu$ error locations plus  zero
counted with multiplicity $t-\mu$; more precisely
  $$\cL(z,{s}_1,\ldots, {s}_{n-k})=(z-z_1)\cdots (z-z_\mu)z^{t-\mu}.$$
 \end{quote}
They labelled such a polynomial {\em general error locator polynomial (GELP)} and proved that it is computable via Gr\"obner computation of the 
syndrome variety ideal ${\Cal I}(\mathbf{V}).$

Actually, they
denoted $\mathbf{V}_{OS}\subset\mathbf{V}$ the set of the non-spurious points of the syndrome variety and considered the polynomial set
$$
{\mathcal{F}_{OS}}=\{f_i,h_j,\chi_i , \lambda_j, {\sf p}_{l\tilde{l}},
\,\,1\leq l< \tilde{l} \leq t, 1 \leq i \leq n-k, 1\leq j\leq t \}\,\,\, \subset
{\cal P},$$ where \vspace{-0.2cm}
\begin{align*}
&f_i :=  \sum\limits_{l=1}^t y_lz_l^{j} - x_{i}, \;\; {\sf p}_{l\tilde{l}}:=  z_{\tilde{l}}  z_l  \frac{z_l^n-z_{\tilde{l}}^n}{z_l-z_{\tilde{l}}}, \\
&h_j :=   z_j^{n+1} -z_j, \;\;\lambda_j :=   y_j^{q-1} -1, \;\;\chi_i :=   x_i^{q^m} -x_i,
\end{align*}
proving that

\begin{Theorem} \cite{locator} It holds $\mathbb{I}({\mathcal{F}_{OS}})= \mathcal{I}(\mathbf{V}_{OS}).$
\end{Theorem}

Moreover they considered the reduced \GR \ basis $G$ of $\mathbb{I}({\mathcal{F}_{OS}})= \mathcal{I}(\mathbf{V}_{OS})$  w.r.t. the lex ordering
with $x_1 < \cdots < x_{n-k} < z_t < \cdots < z_1 < y_1 < \cdots <
y_t$ and, adapting Gianni-Kalkbrenner Theorem \cite{gianni,kalk}, denoted,
 for each $\iota\leq t$ and  each $\ell\in{\Bbb N}$
 $$G_\iota
:= G\cap{\Bbb F}_q[x_1,\ldots,x_{n-k},z_t ,\cdots, z_\iota] \And
G_{\iota\ell} := \{g\in G_\iota\setminus G_{\iota+1} : \deg_{x_\iota}(g) = \ell\};$$
further enumerating each $G_{\iota\ell}$ as
$$G_{\iota\ell} := \{g_{\iota\ell1},\ldots,g_{\iota\ell j_{\iota\ell}}\},
{\bf T}(g_{\iota\ell1})  < \cdots <  {\bf T}(g_{\iota\ell
j_{\iota\ell}}),$$
they proved

\begin{Theorem}\label{OSALA1} \cite{locator} With the present notation we have
\begin{enumerate}
\item $G \cap \FF_{q}[x_1,...,x_{n-k},z_1 , \dots, z_t] = \cup_{i=1}^{t}G_i$;

\item $G_i=  \sqcup_{\delta = 1}^{i}G_{i\delta}$ and $G_{i \delta}
 \neq \emptyset $, $1 \leq i \leq t$, $1 \leq \delta \leq i$;

\item $G_{ii}=\{g_{ii1}\}$, $1 \leq i \leq t$, i.e. exactly one
polynomial exists with degree $i$ w.r.t. the variable $z_i$ in
$G_i$;

\item  ${\bf T}(g_{ii1})=z_i^i,\,\,\,\,\,\, \lc(g_{ii1})=1$;

\item if $ 1 \leq i \leq t $ and $ 1 \leq \delta \leq i-1$, then $\forall
g \in G_{i \delta}, z_i\mid g$.
\end{enumerate}
\end{Theorem}

\begin{Definition} \cite{locator}
The unique polynomial
$$g_{tt1}= z_t^t+\sum_{l=1}^t a_{t-l}({s}_1,\ldots, {s}_{n-k})z_t^{t-l}$$ with degree
$t$ w.r.t. variable $z_t$ in $G_t$, which is labelled the  {\em general error locator polynomial}, is such that
the following properties are equivalent for each syndrome vector ${s}=({s}_1,\ldots, {s}_{n-k}) \in ({\Bbb
F}_{q^m})^{n-k}$ corresponding to an error with weight bounded by $t$:
\begin{itemize}
\item there are exactly $\mu\leq t$ errors $\zeta_1,\ldots\zeta_\mu$;
\item $a_{t-l}({s}_1,\ldots, {s}_{n-k}) = 0$ for $l>\mu$ and
$a_{t-\mu}({s}_1,\ldots, {s}_{n-k}) \neq 0$;
\item $g_{tt1}({{s}_1,\ldots, {s}_{n-k}},z_t) =
z^{t-\mu}\prod_{j=1}^\mu (z-\zeta_i)$.
\end{itemize}
This means
that the general error locator polynomial
$g_{tt1}$ is the monic polynomial in $\FF_{q}[x_1,...,x_{n-k},z]$ which
satisfies the following property:
\begin{itemize}
\item[]  {\em given
a syndrome vector ${s}=({s}_1,\ldots, {s}_{n-k}) \in ({\Bbb
F}_{q^m})^{n-k}$ corresponding to an error with weight $\mu \leq
t$, then its $t$ roots are the $\mu$ error locations plus  zero
counted with multiplicity $t-\mu$.}
\end{itemize}

\end{Definition}
\begin{Theorem}[\cite{locator}]
Every cyclic code possesses a general error locator polynomial.
\end{Theorem}

\subsection{Zech Tableaux}
We observe that the parameters of a minimal basis $$G=\{t_1,\ldots,t_r\}, t_i:=x_1^{a_i}x_2^{b_i}, t_1<\cdots<t_r$$ of a monomial ideal 
$${\sf T}\subset{\Cal T}=\{x_1^{a}x_2^b : (a,b)\in{\Bbb N}^2\}$$
satisfy relations
\begin{itemize}
 \item $a_1>a_2>\cdots>a_r$
 \item $b_1<b_2<\cdots<b_r$
 \item and ${\sf T}$ is 0-dimensinal if and only if $b_1=0=a_r$.
 \end{itemize}
The corresponding escalier ${\sf N}={\Cal T}/{\sf T}$ is
$${\sf N}:=\bigsqcup_{i=1}^{r-1}\{x_1^{a}x_2^b : 0\leq a<a_i,0\leq b<b_{i+1}\}.$$

Let us now consider the field ${\mathbb F}_{2^m}={\mathbb Z}_2[a]$, $a$ denoting a primitive $(2^m-1)^{\rm th}$ root of unity; for a value $n\mid (2^m-1)$ we denote $\alpha:=\frac{2^m-1}{n}$ and 
$b:=a^{\alpha}$ a primitive $n^{\rm th}$ root of unity. 
Denote, for $i, 0\leq i<\alpha$, $Z_i :=\{j, 1\leq j\leq n : 1+b^j=1+a^{j\alpha}\equiv a^{i \bmod n}\}$, set $z(i) =\#Z_i$; for any set $S\subset\{j, 1\leq j\leq n\}$  we consider also the values  $\zeta(i)=\#(S\cap Z_i)$ and we call  {\em $((2^m-1),n;S)$-Zech Sequence} the sequence
$\left(\zeta(0),\zeta(1),\ldots,\zeta(\alpha-1)\right)$.

\begin{Definition}\cite{Nadir} In this contex, let us consider a
set $H\subset\{j, 1\leq 
j\leq n\}$ and the corresponding $((2^m-1),n;H)$-Zech Sequence  
$\left(\zeta(0),\zeta(1),\ldots,\zeta(\alpha-1)\right)$.
\\
 The {\em $(2^m-1,n;H)$-Zech Tableau} is the assignement of
\begin{itemize}
 \item an ordered sequence $S:=[j_0,...,j_{r-1}]\subset\{i, 0\leq i<\alpha\}$ which satisfies
 \begin{itemize}
 \item $\zeta(j_0)\geq\ldots\geq \zeta(j_{r-1})>0$,
 \item $\zeta(j)=0 \Forall j\notin S$.
\end{itemize}
 \item the minimal basis ${\sf G}=\{t_1,\ldots,t_r\}, t_i:=x_1^{a_i}x_2^{b_i},$ of the monomial ideal ${\sf T}={\Cal T}\setminus{\sf N}$ corresponding to the escalier 
 $${\sf N}:=\{x_1^{a}x_2^b : 0\leq a<r,0\leq b<b_{\zeta(j_a)}\}.$$
\end{itemize}
\end{Definition}
\section{The error variety}
Our aim is to describe the syndrome variety of a  binary cyclic $[n,2,d]$-code  
$C$ over $GF(2^m)$, length $n\mid 2^m-1$ and primary defining set $S_C=\{1,l\}$.
Thus we consider
\begin{itemize}
 \item the field ${\mathbb F}_{2^m}={\mathbb Z}_2[a]$, $a$ denoting a primitive $(2^m-1)^{\rm th}$ root of unity, $\alpha:=\frac{2^m-1}{n}$ and 
$b:=a^{\alpha}$ a primitive $n^{\rm th}$ root of unity, 
\item the points
$(c+d,c^l+d^l,c,d)\in\left({\mathbb F}_{2^m}\right)^4 : c:=b^{\ind(c)}, d:=b^{\ind(d)}\in\{e\in{\mathbb F}_{2^m} : e^n=1\}=:{\Cal R}_n,$ and
\item the set ${\Cal S}_n :={\Cal R}_n\sqcup\{0\}$.
\end{itemize}

Let us consider the following data
\begin{description}
\item[] ${\sf Z}_e^\times :=\{(c,d), c,d\in {\Cal S}_n\}$,$\#{\sf Z}_e^\times=(n+1)^2$, the set of all error locations corresponding to errors of weight $w, 0\leq w\leq 2$,
\item[] ${\sf Z}_{ns}^\times :=\{(c,d), c,d\in {\Cal S}_n\}\setminus\{(c,c),c\in{\Cal R}_n\}$,$\#{\sf Z}_{ns}^\times=n^2+n+1$, the set of all {\em non spurious} error locations corresponding to errors of weight $w, 0\leq w\leq 2$,
\item[] ${\sf Z}_+^\times :=\{(c,d), c,d\in{\Cal S}_n, c\neq d\}$,$\#{\sf N}_+^\times=n^2+n$, the set of all non-spurious  error locations  corresponding to errors of weight $w, 1\leq w\leq 2$,
\item[] ${\sf Z}_2^\times :=\{(c,d), c,d\in{\Cal R}_n, c\neq d\}$,$\#{\sf N}_2^\times=n^2-n$, the set of all non-spurious  error locations  corresponding to errors of weight $w= 2$;
\item[] ${\sf N}_e^\times :=\{z_1^{i}z_2^{j}, 0\leq i\leq n, 0\leq j\leq n\}$,$\#{\sf N}_e^\times=(n+1)^2$,
\item[] ${\sf N}_{ns}^\times :=\{z_1^{i}z_2^{j}, 0\leq i\leq n, 0\leq j< n,\}\sqcup\{z_2^{n}\}$,$\#{\sf N}_{ns}^\times=n^2+n+1$,
\item[] ${\sf N}_+^\times :=\{z_1^{i}z_2^{j}, 0\leq i\leq n, 0\leq j< n\}$,$\#{\sf N}_+^\times=n^2+n$,
\item[] ${\sf N}_2^\times :=\{z_1^{i}z_2^{j}, 0\leq i< n, 0\leq j< n-1\}$,$\#{\sf N}_e^\times=n^2-n$;

\item[] ${G}_{e}^\times :=\{z_1^{n+1}-z_1,z_2^{n+1}-z_2 \}$,
\item[] ${G}_{ns}^\times :=\{z_1^{n+1}-z_1,z_1\left((z_2-z_1)^{n}-1\right),z_2^{n+1}-z_2\}$,
\item[] ${G}_{+}^\times :=\{z_1^{n+1}-z_1,(z_2-z_1)^{n}-1\}$, 
\item[] ${G}_2^\times :=\{z_1^{n}-1,z_2^{n-1}-\sum_{i=1}^{n-1} \binom{n-1}{i}z_1^iz_2^{n-1-i}\}$.
\end{description}

For $\ast\in\{e,ns,+,2\}$ setting ${\sf I}_\ast^\times :={\Bbb I}(G_\ast^\times)$ the ideal generated by $G$, as a direct corollary of Macaualay's trick and Lazard correspondence,
we trivially have the following facts:
\begin{Corollary} It holds

\begin{enumerate}
\renewcommand\theenumi{{\rm (\alph{enumi})}} 
\item ${\sf N}_\ast^\times$ is the Gr\"obner escalier of ${\sf I}_\ast^\times$;
\item for each degree compatible term-ordering s.t. $z_2>z_1$, $G_\ast^\times$ is the Gr\"obner basis of the ideal ${\sf I}_\ast^\times$ it generates;
\item ${\sf Z}_\ast^\times={\Cal V}({\sf I}_\ast^\times)$ and 
${\sf I}_\ast^\times={\Frak I}({\sf Z}_\ast^\times)$.
\end{enumerate}
\end{Corollary}

\section{The syndrome variety}

Our aim being decoding a  binary BCH $[n,2,d]$-code  
$C$ over $GF(2^m)$, length $n\mid 2^m-1$ and primary defining set $S_C=\{1,l\}$, we begin by reformulaing in our setting some preliminary notations.

In this setting the syndrome variety would be specialized to the set
$$\mathbf{V}_2:=\left\{(c+d,c^l+d^l,c,d)\in {\Cal S}_n^{4}\right\},$$
the syndrome map to the map 
$$\phi :   {\Cal S}_n^2\to  {\Cal S}_n^2 : (c,d)\mapsto (c+d,c^l+d^l),$$
and GELP becomes the polynomial
$$\cL(x_1,x_{2},z_1)=z_1^2+\sum_{i=1}^2a_{2-i}(x_1,x_{2})z_1^{2-i} : \cL((c+d,c^l+d^l,z_1)=(z_1-c)(z_1-d)\in \FF_q[x_1,x_{2}][z_1].$$
and could be obtained, no more by Gr\"obner basis computation, but via Lagrange interpolation, obtainable by M\"oller algorithm \cite{MMM2} or by Lundqvist interpolation formula \cite{Lu} of the polynomial
$\sum\limits_{t\in{\sf N}({\Frak I}(\mathbf{V}_2))} c_tt$ over the set $\mathbf{V}_2$, {\em provided that we have ${\sf N}({\Frak I}(\mathbf{V}_2))$.}

This is the first motivation of this note.
Considering 
the following sets of points
\begin{description}
\item[] ${\sf Z}_2 :=\{(c+d,c^l+d^l,c,d), c,d\in{\Cal R}_n, c\neq d\}$,$\#{\sf Z}_2^\times=n^2-n$;
\item[] ${\sf Z}_+ :=\{(c+d,c^l+d^l,c,d), c,d\in{\Cal S}_n, c\neq d\}$,$\#{\sf Z}_+^\times=n^2+n$,
\item[] ${\sf Z}_{ns} :=\{(c+d,c^l+d^l,c,d), c,d\in {\Cal S}_n\}\setminus\{(0,0,c,c),c\in{\Cal R}_n\}$,$\#{\sf Z}_{ns}^\times=n^2+n+1$,
\item[] ${\sf Z}_e :=\{(c+d,c^l+d^l,c,d), c,d\in {\Cal S}_n\}$,$\#{\sf Z}_e^\times=(n+1)^2$,
\end{description}
we show that, assuming known ${\sf N}_2$ direct application of Lazard's Cerlienco-Mureddu Correspondence and some trivial considerations on the syndrome map $\phi$ allow to deduce ${\sf N}_\ast$, $\Phi_\ast$ and also the irrelevant ${G}_\ast$ for each h $\ast\in\{e,ns,+,2\}.$

The shape of ${\sf N}_2$ for a $[2^m-1,2]$-code $C$ over $\FF_{2^m}$ has been given recently by 
\cite{AcaHelp} which moreover
proved that
the related \GR\ basis has the shape
$$G=(x_1^n-1,g_2,z_2+z_1+x_1,g_4)$$
where (see \cite{locator})
$$g_2=\frac{x_2^{\frac{n+1}{2}}-x_1^\frac{n+1}{2}}{x_2-x_1}=
x_2^{\frac{n-1}{2}}+\sum_{i=1}^{\frac{n-1}{2}} \binom{\frac{n-1}{2}}{i}x_1^ix_2^{\frac{n-1}{2}-i}$$
and $g_4=z_1^2-\sum_{t\in{\bf N}} c_tt$ is Sala-Orsini general error locator polynomial.

On the basis of that, an improvement of  GELP was introduced, under the label of {\em half error locator polynomial}  (HELP), in \cite{AcaHelp}, which remarks that the relation
$x_1+z_1+z_2\in G_\ast$, requires just to produce a polynomial
$${\Cal H}(x_1,x_2)=\sum_{t\in{\sf N}^x_e} c_tt : {\Cal H}(c+d,c^l+d^l)=c \Forall (c,d)\in{\Cal R}_n^2.$$
Such polynomial could be obtained by fixing randomly an element $c$ in each pair $(c,d)\in{\Cal R}_n^2$ and apllying Lagrange/Lundqvist interpolation on the set $\left\{(c+d,c^l+d^l,c): (c,d)\in{\sf Z}_2\right\}$.

A more efficient solution has been proposed in \cite{Help} which suggests to reformulate 
${\sf Z}_2$  as
$${\sf Z}_2 :=\{(c(1+b^{2i-1},c^l(1+b^{l(2i-1)},c,cb^{2i-1}), c\in {\Cal R}_n,\}, 0\leq i<\frac{n-1}{2}.$$

\subsection{Trivalities on the syndrom map}

We have
\begin{enumerate}
\item $\phi(c,d)=\phi(d,c)$;
\item $\#\Jm(\phi)= \frac{(n+1)^2}{2}$;
\item $\phi(c,d) = (0,0)$ if and only if either $(c,d)=(0,0)$ {\em id est} $w=0$ or $0\neq c = d$ id est $(c,d)$ is spurious;
\item $\phi(c,d)= (x, x^l),\, x\neq0,$ if  and only if  $cd\sum_{i=0}^{l-2} c^id^{l-i-2}=0$ .
\begin{proof}
Since $\phi(c,d)=(c+d,c^l+d^l)$ 
$$\phi(c,d)= (x, x^l) \iff c^l+d^l=(c+d)^l=\sum_{i=0}^{l} c^id^{l-i-2}\iff cd\sum_{i=0}^{l-2} c^id^{l-i-2}=0$$
\end{proof}

\item In particular, in case $l=3$, $\phi(c,d)= (x, x^3),\, x\neq0,$ if and only if $cd=0$ if and only if $0<w<2.$

\item for each $x\in GF(2^m)^*, (0,x)\notin\Jm(\phi)$, since 
$$c+d=0\iff c=d \iff (c,d) \mbox{ is sporious } \then c^l=d^l.$$

\item for each $x\in GF(2^m)^*, (x,0)\in\Jm(\phi)$ then $\gcd(n,l)\neq 1$
\begin{proof}
Since $\phi(c,0)=(c,c^l)$ and $c^l=0 \then c=0$, we are restricted to consider pairs
$(c,d)$ with $c\neq0\neq d \neq c$.
Then  $t:= \frac{d}{c}\in GF(2^m)^*$ and
$d^l=c^l\iff t^l=1\iff \ord(t)\mid l$; since $n=\ord(t)\ind(t)$, $ \ord(t)\mid l\iff\ord(t)\mid\gcd(n,l)$. Thus  
\end{proof}

\end{enumerate}

\subsection{$n=2^m-1$} Let us now fix 
the lex ordering with  $x_1<x_2<z_1<z_2$ and let us consider  
the ideal ${\sf I}_\ast:={\Frak I}({\sf Z}_\ast)$ 
and remark that, if $n=2^m-1$, so that $\alpha=1$ since \cite{AcaHelp}
\begin{enumerate}
\renewcommand\theenumi{{\rm (\alph{enumi})}}  
\item ${\sf N}_2 :=\{z_1^kx_1^{i}x_2^{j}, 0\leq i< n, 0\leq j< \frac{n-1}{2},0\leq k\leq1\}$,$\#{\sf N}_2=n^2-n$,
\end{enumerate}
we can deduce that
\begin{enumerate}
\renewcommand\theenumi{{\rm (\alph{enumi})}} 
\setcounter{enumi}{1}
\item ${\sf N}_2$ corresponds with ${\sf Z}_2$ via the Cerlienco Mureddu Correspondence
 $$\Phi_2(x_1^{i}x_2^{j})=(b^{i}+b^{j},b^{li}+b^{lj},b^{i},b^{j}), 
 \Phi_2 (z_1x_1^{i}x_2^{j})=(b^{i}+b^{j},b^{li}+b^{lj},b^{j},b^{i}),$$
\item ${G}_2 :=\{x_1^{n}-1,x_2^{\frac{n-1}{2}}-\sum_{i=1}^{\frac{n-1}{2}} \binom{\frac{n-1}{2}}{i}x_1^ix_2^{\frac{n-1}{2}-i},z_2+x_1+x_2,z_1^2- \sum_{t\in{\sf N}_2} c_tt\}$;
\item ${\sf Z}_+ = {\sf Z}_2\sqcup\{(c,c^l,c,0),(c,c^l,0,c), c\in{\Cal R}_n\}$;
\item ${\sf N}_+ = {\sf N}_2\sqcup\{z_1^kx_1^{i}x_2^{\frac{n+1}{2}}, 0\leq i< n, 0\leq k\leq1\}$,
\item corresponding to the Cerlienco Mureddu Correspondence 
$$\Phi_+(x_1^{i})=(b^{i},b^{li},b^{i},0), 
\Phi_+ (z_1x_1^{i}) =(b^{i},b^{li},0,b^{i}) 
\And \Phi_+(z_1^kx_1^{i}x_2^{j})=\Phi_2(z_1^kx_1^{i}x_2^{j+1}) \Forall z_1^kx_1^{i}x_2^{j}\in{\sf Z}_2$$
\item ${G}_+ :=\{x_1^{n}-1,x_2\left(x_2^{\frac{n-1}{2}}-\sum_{i=1}^{\frac{n-1}{2}} \binom{\frac{n-1}{2}}{i}x_1^ix_2^{\frac{n-1}{2}-i}\right),z_2+x_1+x_2,z_1^2- \sum_{t\in{\sf N}_+} c_tt\}$,
\item ${\sf Z}_{ns} = {\sf Z}_+\sqcup\{(0,0,c,c), c\in{\Cal R}_n\}$;
\item ${\sf N}_{ns} = {\sf N}_+\sqcup\{x_1^{n},z_1x_1^{n}\}\sqcup\{z_1^{k}, 2\leq k<n\}$,
\item corresponding to the Cerlienco Mureddu Correspondence 
\begin{eqnarray*}
 \Phi_{ns}(z_1^{k})=(0,0,b^{k},b^{k}),&& 0\leq k<n, \\
 \Phi_{ns}(z_1^kx_1^{i})=\Phi_+(z_1^kx_1^{i+1}), && \Forall z_1^hx_1^{i}\in{\sf Z}_2, \\
 \Phi_{ns}(z_1^kx_1^{i}x_2^{j})=\Phi_+(z_1^kx_1^{i}x_2^{j}), && \Forall z_1^hx_1^{i}x_2^{j}\in{\sf Z}_2, j\neq0.\\
\end{eqnarray*}
\item ${G}_{ns} :=\{x_1^{n+1}-x_1,x_2\left(x_2^{\frac{n-1}{2}}-\sum_{i=1}^{\frac{n-1}{2}} \binom{\frac{n-1}{2}}{i}x_1^ix_2^{\frac{n-1}{2}-i}\right),z_2+x_1+x_2,
x_1\left(z_1^2- \sum_{t\in{\sf N}_+} c_tt\right),z_1^n-1\}$,
\item ${\sf Z}_e = {\sf Z}_{ns}\sqcup\{(0,0,0,0),\}$;
\item ${\sf N}_e = {\sf N}_{ns}\sqcup\{z_1^{n}\}$,
\item corresponding to the Cerlienco Mureddu Correspondence 
$$\Phi_{e}(z_1^{n})=(0,0,0,0) 
\And \Phi_{e}(x_2^kz_1^{i}z_2^{j})=\Phi_{ns}(x_2^kz_1^{i}z_2^{j}) \Forall x_2^hz_1^{i}z_2^{j}\in{\sf Z}_{ns}$$
\item ${G}_{e} :=\{x_1^{n+1}-x_1,x_2\left(x_2^{\frac{n-1}{2}}-\sum_{i=1}^{\frac{n-1}{2}} \binom{\frac{n-1}{2}}{i}x_1^ix_2^{\frac{n-1}{2}-i}\right),z_2+x_1+x_2,
x_1\left(z_1^2- \sum_{t\in{\sf N}_+} c_tt\right),z_1^{n+1}-z_1\}$.
\end{enumerate}

\subsection{$n\neq2^m-1$}\label{XAX} If we now assume that $n\neq2^m-1$, so that $\alpha>1$, the argument developped above can be {\em verbatim} reformulated provided that we know
\begin{enumerate}
\renewcommand\theenumi{{\rm (\alph{enumi})}}  
\item the structure of the order ideal ${\sf N}_2$,$\#{\sf N}_2=n^2-n$, i.e.  
a minimal basis $\{t_1,\ldots,t_r\}, t_i:=x_1^{a_i}x_2^{b_i},$ of the monomial ideal, 
\item a Cerlienco Mureddu Correspondence $\Phi_2: {\sf N}_2 \to {\sf Z}_2$,
\item the Gr\"obner basis $G_2 =\{g_1,\ldots,g_r\}, {\bf T}(g_i) = t_i:=x_1^{a_i}x_2^{b_i}, t_1<\cdots<t_r$.
\end{enumerate}
We in fact obtain
\begin{enumerate}
\renewcommand\theenumi{{\rm (\alph{enumi})}} 
\setcounter{enumi}{3}
\item ${\sf Z}_+ = {\sf Z}_2\sqcup\{(c,c^l,c,0),(c,c^l,0,c), c\in{\Cal R}_n\}$;
\item ${\sf N}_+ = \{x_2t, t\in{\sf N}_2\}\sqcup\{z_1^kx_1^{i}, 0\leq i< n, 0\leq k\leq1\}$,
\item corresponding to the Cerlienco Mureddu Correspondence 
$$\Phi_+(x_1^{i})=(b^{i},b^{li},b^{i},0), 
\Phi_+ (z_1x_1^{i}) =(b^{i},b^{li},0,b^{i}) 
\And \Phi_+(tx_2)=\Phi_2(t) \Forall t\in{\sf Z}_2;$$
\item ${G}_+ :=\{g_1,x_2g_2,\ldots,x_2g_r,z_2+x_1+x_2,z_1^2- \sum_{t\in{\sf N}_+} c_tt\}$,
\item ${\sf Z}_{ns} = {\sf Z}_+\sqcup\{(0,0,c,c), c\in{\Cal R}_n\}$;
\item ${\sf N}_{ns} = {\sf N}_+\sqcup\{x_1^{a_1},z_1x_1^{a_1}\}\sqcup\{z_1^{k}, 2\leq k<n\}$,
\item corresponding to the Cerlienco Mureddu Correspondence 
\begin{eqnarray*}
 \Phi_{ns}(z_1^{k})=(0,0,b^{k},b^{k}),&& 0\leq k<n, \\
 \Phi_{ns}(z_1^kx_1^{i})=\Phi_+(z_1^kx_1^{i+1}), && \Forall z_1^hx_1^{i}\in{\sf Z}_2, \\
 \Phi_{ns}(z_1^kx_1^{i}x_2^{j})=\Phi_+(z_1^kx_1^{i}x_2^{j}), && \Forall z_1^hx_1^{i}x_2^{j}\in{\sf Z}_2, j\neq0.\\
\end{eqnarray*}
\item ${G}_{ns} :=\{x_1g_1,x_2g_2,\ldots,x_2g_r,z_2+x_1+x_2,
x_1\left(z_1^2- \sum_{t\in{\sf N}_+} c_tt\right),z_1^n-1\}$,
\item ${\sf Z}_e = {\sf Z}_{ns}\sqcup\{(0,0,0,0),\}$;
\item ${\sf N}_e = {\sf N}_{ns}\sqcup\{z_1^{n}\}$,
\item corresponding to the Cerlienco Mureddu Correspondence 
$$\Phi_{e}(z_1^{n})=(0,0,0,0) 
\And \Phi_{e}(t)=\Phi_{ns}(t) \Forall t\in{\sf Z}_{ns}$$
\item ${G}_{e} :=\{x_1g_1,x_2g_2,\ldots,x_2g_r,z_2+x_1+x_2,
x_1\left(z_1^2- \sum_{t\in{\sf N}_+} c_tt\right),z_1^{n+1}-z_1\}$.
\end{enumerate}

\begin{Remark} \noindent
\begin{enumerate}
\item In 4.3.(a) necessarily $b_1=a_r=0$ which is sufficient to justify the formulas
 ${G}_+\cap\FF_q[x_1,x_2] =\{g_1,x_2g_2,\ldots,x_2g_r\}$ in 4.3.(g) and ${G}_{ns}\cap\FF_q[x_1] =\{x_1g_1\}.$
 \item In  4.2.(e) we needed to add to ${\sf N}_2$ $2n$ new elements and a natural application of Cerlienco-Mureddu Algorithm introduced the monomial $\{z_1^kx_1^{i}x_2^{\frac{n+1}{2}}, 0\leq i< n, 0\leq k\leq1\}$ but in that case ${\sf N}_2$ had the shape of two superimposed rectangles and it was sufficient to add two further strips of monomial over the monomial $\{z_1^kx_1^{i}x_2^{\frac{n-1}{2}}$ as explicitly stated in 4.2.(f). 
In the case 4.3.(e) the shape of the two superimposed terms is not available and the only natural solution is to apply Cerlienco-Mureddu Algorithm considering first the $2n$ new elements and later the elements of ${\sf N}_2$ and this fact too is explicitly recordered in 4.3.(f).
\end{enumerate}
 
\end{Remark}

\end{document}